# BLOW-ANALYTIC RETRACTION ONTO THE CENTRAL FIBRE

Adam Parusiński

ABSTRACT. Let $\mathcal{X}$ be a complex analytic space and let $f: \mathcal{X} \to \mathbf{C}$ be a proper complex analytic function with nonsingular generic fibres. By adapting the blow-analytic methods of Kuo we construct a retraction of a neighbourhood of the central fibre $f^{-1}(0)$ onto $f^{-1}(0)$. Our retraction is defined by the flow of a real analytic vector field on an oriented real analytic blow-up of $\mathcal{X}$. Then we describe in terms of this blow-up the associated specialization map and local Milnor fibrations. The method also works in real analytic category.

## 0. Introduction

Let $\mathcal{X}_0$ be a compact complex analytic space and let $f: \mathcal{X} \to \mathbf{C}$ be a one parameter deformation of $\mathcal{X}_0$, that is $f^{-1}(0) = \mathcal{X}_0$. We call $\mathcal{X}_0$ the central fibre of $f$. The pair $\mathcal{X}, \mathcal{X}_0$ can be triangulated, and hence $\mathcal{X}_0$ is a strong deformation retract of its neighbourhood in $\mathcal{X}$. Thus shrinking $\mathcal{X}$, if necessary, we may assume that there is a retraction $r: \mathcal{X} \to \mathcal{X}_0$. Often, one would like to have retractions satisfying special properties, for instance such that the restriction of $r$ to a general fibre $\mathcal{X}_s = f^{-1}(s)$, $s \ne 0$ and small, called a specializaton map $r_s: \mathcal{X}_s \to \mathcal{X}_0$, does not depend on the choice of $s$, and that the fibres of the specialization map are constant over the strata of some stratification of $\mathcal{X}_0$. The purpose of this paper is to present a method of construction of such retraction under an additional assumption that $f$ is a smoothing, that is the general fibres of $f$ are nonsingular. Our method works also for real analytic functions.

In his exposé on vanishing cycles [De] Deligne assumes the existence of such good retraction stating that "il me semble plausible" that it exists. Actually, to construct the sheaf of vanishing (or rather neighbouring cycles) he uses specialization defined up to a homotopy, that is $\mathcal{X}_s \hookrightarrow \mathcal{X} \hookleftarrow \mathcal{X}_0$, which works since the latter inclusion is a homotopy equivalence.

1991 *Mathematics Subject Classification.* 32C05, 32C15, 32S30, 14P15.
Dedicated to Professor Tzee Char Kuo for his sixtieth birthday
The author would like to thank Roy Smith and Robert Varley for suggesting the problem, and also Clint McCrory for many helpfull discussions during the preparation of this paper

Typeset by $\mathcal{A}_{\mathcal{M}}\mathcal{S}$-TeX

One may find in the literature some examples of actual constructions of retractions which use diverse techniques. As we mentioned above triangulability of analytic sets suggests one possible approach. Note that in order to obtain a resonable description of the specialization map one has to triangulate $f$ as a mapping. This was carried out by Shiota in the real case in [Sh1], in the complex algebraic one in [Sh2, §3], and in the local complex analytic one in [Sh2, §2]. A construction of similar type was obtained by semi-algebraic methods in [CR]. Note that PL or semi-algebraic methods ensure the existence of semi-algebraic (subanalytic) stratification of the central fibre $\mathcal{X}_0$, such that both the retraction and the specialization map has constant type on the strata of such stratification. T his seems to be not sufficient for complex analytic or algebraic functions. Another possible approach, working in both complex and real case, is the use of Whitney and Thom stratifications. This idea was carried out in [GM, §6] by the use of control functions and controlled vector fields. The retraction constructed in [GM] is actually not a retraction in the strict meaning because its restriction to $\mathcal{X}_0$ is only homotopic to identity.

Note that triangulation and stratification are the most common methods of studying and establishing topological triviality of families of real algebraic (analytic) functions. There is another interesting equisingularity condition called the blow-analytic equivalence. It originated in works of Kuo (see, for instance [Kuo1], [Kuo2]) and is based on resolution of singularities and construction of real analytic vector field on the blow-up space. The trivializations obtained by this method, though not analytic, blow up to real analytic homeomorphisms, and hence are called blow-analytic. In this paper we use Kuo's approach to construct a retraction satisfying the desired properties.

Let $f\colon \mathcal{X} \to \mathbf{R}$ be a real analytic function defined on a compact real analytic space $\mathcal{X}$. As before we shall assume that the generic fibres of $f$ are smooth. Then we may find $\sigma\colon M \to \mathcal{X}$ a modification of $\mathcal{X}$, such that $M$ is nonsingular, $\sigma$ is an isomorphism over $\mathcal{X} \setminus \mathcal{X}_0$, and $f \circ \sigma$ is normal crossings. Then, we construct a retraction of $M$ onto $X = \sigma^{-1}(\mathcal{X}_0)$, which is easier since $X$ is a divisor with normal crossings. This we achieve first 'cutting $M$ along $X$' and then retract along the trajectories of a real analytic vector field. Finally we push the retraction of $M \to X$ down to $\mathcal{X}$ to obtain a desired retraction $\mathcal{X} \to \mathcal{X}_0$. We describe this construction in details in Section 3.

In the complex case, we again resolve singularities of $\mathcal{X}$ and make the central fibre a divisor with normal crossings. Then, it is easy to give a formula for retraction in terms of local coordinates. Thus the problem reduces to patching such local retractions together. This was achieved in [Cl] by a partition of unity. Instead, we use the idea of A'Campo [A'C] and blow up the complex analytic divisor as a real analytic subspace. Then we apply the real analytic case. Note that, such real analytic blowing up followed by our 'cutting along' real analytic divisor can be identified with the oriented blowing up of [A'C]. The complex case is described in Section 4.

Obviously, one would like to get rid of the assumption that $f$ is a smoothing, that is that the general fibre of $f$ is nonsingular. To extend our construction to the general case one has to overcome similar difficulties as these which cause that Kuo's equisingularity criterion is established only for families of functions with isolated singularities, see [Kuo2], Conjecture of §2 and the last conjecture of §3, in particular.

## 1. Some properties of real analytic (sub)manifolds

In this section we recall some basic properties of real analytic hypersurfaces of real analytic manifolds. These properties can be derived easily from the embeddingtheorem of Grauert [Gr] and from the Theorems A and B of Cartan [Car]. We restrict ourselves to the case of nonsingular hypersurfaces even if most of the facts stated below hold also for singular coherent subspaces of codimension 1.

Let $M$ be a paracompact real analytic manifold of pure dimension $n$. By [Gr] we can always assume that $M$ is a closed analytic submanifold of $\mathbf{R}^N$ for $N$ sufficiently big. Each nonsingular hypersurface $V$ of $M$ defines a coherent sheaf of analytic functions vanishing on $V$ which can be identified with the sheaf of sections of an analytic line bundle $L(V)$. In particular, $V$ is *the zero subspace* of a real analytic section $\mathbf{e}$ of $L(V)$.

**Proposition 1.1.** *Let $M$ be a paracompact real analytic manifold and let $V$ be a nonsingular hypersurface of $M$. Then:*
(a) *$V$ is the zero subspace of an analytic section $\mathbf{e}$ of an analytic line bundle $L = L(V)$ over $M$. The fundamental class $[V] \in \mathbf{H}_1^{cl}(M; \mathbf{Z}_2)$ is dual to the first Stiefel-Whitney class $w_1(L) \in \mathbf{H}^1(M; \mathbf{Z}_2)$ of $L$.*
(b) *$L$ admits an analytic metric, so $V$ as a set is defined by an analytic function $|\mathbf{e}|^2$.*
(c) *The following conditions are equivalent:*
   (1) *$L$ is topologically trivial;*
   (2) *$L$ is trivial as an analytic line bundle;*
   (3) *$w_1(L) = 0$;*
   (4) *$V$ is the zero **subspace** of an analytic function;*
(d) *if $V$ is connected then the conditions of (c) are equivalent to*
   (5) *$V$ disconnects the connected component of $M$ which contains $V$.*

*Proof (sketch).* (a) is standard. (b) follows from [Car]. To prove (c) we note that: (1) $\iff$ (2) follows from [Car, 7.7]. The rest of (c) is easy Assume that $V$ is connected. We shall show (4) $\iff$ (5). Let $M_1$ be a connected component of $M$ containing $V$. If $V$ is defined as a subspace by an analytic function $f$, then $V$ disconnects $M_1$ into the part where $f$ is negative and that where it is positive. On the other hand, if $V$ disconnects $M_1$, then $[V]$ is a boundary so $w_1(L) = 0$. $\square$

*Remark 1.2.* $L^2 = L \otimes L$ is always trivial.

*Remark 1.3.* If $V$ is the exceptional divisor of a blowing-up (with a nonsingular centre of codimension $\geq 2$), then $L$ is not trivial (in any neighbourhood of $V$) since $M_1 \setminus V$ is (even near $V$) connected.

## 2. Stratifications with normal crossings

Let $M$ be a real analytic manifold and let $X$ be a subset of $M$. We shall assume that $X$ admits some kind of analytic structure by being analytic, semianalytic or in full generality subanalytic. By a *(subanalytic) stratification of $X$* we mean a locally finite collection $\mathcal{S} = \{S_\alpha\}$ of disjoint subanalytic subsets of $M$ such that

$$X = \bigcup S_\alpha,$$

and such that the sets $S_\alpha \in \mathcal{S}$, called *strata*, are locally closed analytic submanifolds of $M$. The strata may be connected or not but we always assume that they are pure dimensional. In this paper we consider mainly a much simpler class of stratifications that is the stratifications with normal crossings which we define below.

**Definition 2.1.** A subset $Q$ of $\mathbf{R}^n$ is a *quadrant* if there is partition of $\{1, \ldots, n\}$ into disjoint subsets $I_0$, $I_+$ and $I_-$ such that

$$Q = \{(x_1, \ldots, x_n) \in \mathbf{R}^n; x_i = 0 \text{ if } i \in I_0, x_j > 0 \text{ if } j \in I_+ \text{ and } x_k < 0 \text{ if } k \in I_-\}.$$

We say that a subset $X$ of real analytic manifold $M$ is a *local quadrant* if each point $p \in M$ admits a coordinate neighbourhood $U$ such that $U \cap X$ is a quadrant.

**Definition 2.2.** We say that a stratification $\mathcal{S}$ of a subset $X$ of a real analytic manifold $M$ is *a stratification with normal crossings* if each stratum of $\mathcal{S}$ is a local quadrant.

The set of singular points of a quadrant (resp. local quadrant) is again a quadrant (resp. local quadrant). Therefore each local quadrant $X$ has a canonical stratification which is with normal crossings. Clearly the local quadrants are the only sets which can be stratified by stratifications with normal crossings.

Among stratifications with normal crossings we introduce an equivalence relation which we call *the piecewise analytic isomorphism* and which is weaker than that induced by analytic isomorphisms. Roughly speaking, we say that the given stratifications of local quadrants $X$, $Y$ are equivalent if there exists a homeomorphism $\Phi\colon X \to Y$ which sends strata onto strata and which restricted to **the closure** of each stratum is an analytic isomorphism. So $\Phi$ can be understood as a collection of analytic isomorphisms defined on the closures of strata such that they vary continuously (and subanalytically) but not necessarily analytically from one stratum to another. We give below a precise definition.

In general, for any analytic manifolds $M, N$ and any $X \subset M$ we call a map $\Phi\colon X \to N$ *analytic* if near each point $p \in X$, $\Phi$ is a restriction of an analytic map defined on a neighbourhood of $p$.

**Definition 2.3.** Let $M, N$ be analytic manifolds and let $\mathcal{S}_X, \mathcal{S}_Y$ be stratifications with normal crossings of $X \subset M$ and $Y \subset N$ respectively. We say that a map $\Phi\colon X \to N$ is *piecewise analytic* if each point $p \in X$ admits a neighbourhood $U$ such that for each stratum $S \in \mathcal{S}_X$ and each connected component $S'$ of $S \cap U$, $\Phi$ restricted to the closure (in $X$) of $S'$ is analytic. We call $\mathcal{S}_X, \mathcal{S}_Y$ *piecewise analytically isomorphic* if there is a homeomorhism $\Phi\colon X \to Y$, sending each stratum of $\mathcal{S}_X$ onto a stratum of $\mathcal{S}_Y$ and such that both $\Phi$ and $\Phi^{-1}$ restricted to strata are piecewise analytic.

For instance, $X = \{(x,y) \in \mathbf{R}^2; xy = 0, x \geq 0 \text{ and } y \geq 0\}$ with canonical stratification and the real line $Y = \mathbf{R}$ are piecewise isomorphic if we choose the stratification of $Y$ as $\{\{t \in \mathbf{R}; t > 0\}, \{t \in \mathbf{R}; t < 0\}, \{0\}\}$.

We note also that each piecewise analytic map is subanalytic (that is its graph is subanalytic) but to be piecewise analytically isomorphic is a much stronger relation than that given by subanalytic homeomorphisms.

## 3. Construction of retraction in real analytic case

Let $M$ be as in Section 1 and let $f\colon M \to \mathbf{R}$ be a proper analytic function. Replacing, if necessary, $f$ by its restriction to $f^{-1}(-\varepsilon, \varepsilon)$, with $\varepsilon > 0$ and sufficiently small, we assume that $0$ is the only possible critical value of $f$. Let $X$ denote the zero set of $f$.

**Definition 3.1.** We say that $f$ is *locally normal crossings* if each point of $M$ admits a coordinate neighbourhood $U$, with coordinates $x = (x_1, ..., x_n)$, such

that

(3.1) $$f(x) = x_1^{a_1} x_2^{a_2} \ldots x_n^{a_n} g(x),$$

where $x \in U$, $a_i \in \mathbf{N} \cup \{0\}$, $g$ is analytic and nowhere vanishes in $U$.

We say that $f$ is *globally normal crossings* if it is locally normal crossings and the zero set $X$ of $f$ is the union of nonsingular hypersurfaces $V_i$, $i = 1, ..., k$, intersecting transversally.

Assume that $f$ is globally normal crossings and $X$ is the union of nonsingular connected hypersurfaces $V_i$, $i = 1, \ldots, k$. Let $L_i, \mathbf{e}_i$ denote the line bundles and their analytic sections corresponding to $V_i$ (by (a) of Proposition 1.1). To construct a nice retraction $r: M \to X$ we consider first the local case.

**3.1. Local case.** Assume that $f$ is given by (3.1) in a ball $U = B(0, R) \subseteq \mathbf{R}^n$ and assume for simplicity that $g$ is positive. We order the coordinate functions in such a way that $a_i > 0$ for $i = 1, \ldots, k$ and $a_j = 0$ for $j > k$. For each $k$-tuple of signs $(\varepsilon_1, \ldots, \varepsilon_k)$, $\varepsilon_i = \pm 1$, we define the retraction $r_{(\varepsilon_1, \ldots, \varepsilon_k)}$ of the quadrant $U(\varepsilon_1, \ldots, \varepsilon_k) = \{x \in U; \varepsilon_i x_i \geq 0, i = 1, \ldots, k\}$ on $X = \{x \in U | f(x) = 0\}$ by projecting along the direction given by the vector

(3.2) $$\mathbf{w}(\varepsilon_1, \ldots, \varepsilon_k) = -(\varepsilon_1 \partial/\partial x_1 + \ldots + \varepsilon_k \partial/\partial x_k).$$

Although this projection is not a restriction of an analytic map we can make it piecewise analytic (in the sense of Section 2) by substratifying $U(\varepsilon_1, \ldots, \varepsilon_k)$ into the union of new quadrants (with respect to a new system of coordinates). We explain it in details for $\varepsilon_1 = \ldots = \varepsilon_k = 1$.

Let $U_1 = U(1, \ldots, 1)$ and let $U_{1,i} = \{x \in U_1; x_i \leq x_j \text{ for } j \in \{1, \ldots, k\} \setminus \{i\}\}$. Clearly $U_{1,i}$ is a quadrant. The $U_{1,i}$'s induce a stratification of $U_1$ with normal crossings. The retraction $r_{(1,\ldots,1)}$ restricted to $U_{1,i}$ is given by

$$x \longrightarrow (x_1 - x_i, \ldots, x_k - x_i, x_{k+1}, \ldots, x_n),$$

and therefore is the restriction of an analytic map. We leave the proof of the following properties of $r_{(1,\ldots,1)}$ to the reader:
  (a) the nonsingular fibres of $f$ are transverse to the induced stratification of $U_1$. Hence this stratification induces a stratification with normal crossings on each fibre $f^{-1}(c)$, $c > 0$;
  (b) let $V_1 = \{x \in U_1; f(x) = 0\}$. Then, for small $c > 0$, $r_{(1,\ldots,1)}$ restricted to $f^{-1}(c)$ is a piecewise analytic isomorphism between $f^{-1}(c)$ and a neighbourhood of the origin in $V_1$.

The retractions $r_{(\varepsilon_1, \ldots, \varepsilon_k)}$ agree on the intersections of their domains and define a retraction $r$ of $U$ on $X$.

**3.2. Cutting $M$ along $X$.** In the global case we construct a retraction in a similar way. First we cut $M$ along $X$ in order to decompose it into a union of local quadrants. Next, on each such quadrant, we define an analytic nonsingular vector field **w** whose flow gives a retraction satisfying the required properties. For the exposition purpose we concentrate on the case when $f$ is everywhere nonnegative leaving the details of the general case to the reader. We start with formal definition of cutting $M$ along $X$.

Let $M$ be as before and let $V$ be a nonsingular analytic hypersurface of $M$. Let $L$ be an analytic line bundle corresponding to $V$ and let **e** be an analytic section of $L$ defining $V$ (as a subspace). By $S(L)$ we denote the space of oriented directions in the fibres of $L$. Then $S(L)$ has a unique structure of real analytic manifold such that the projection $\pi_L \colon S(L) \to M$ is a double analytic covering (trivial if and only if $L$ is trivial).

If we fix an analytic metric on $L$, then $S(L)$ is isomorphic to the sphere bundle of $L$. The induced line bundle $(\pi_L)^*L$ on $S(L)$ is always trivial and $\pi_L^{-1}(V)$ is the zero space of an analytic function $S(L) \ni (p,l) \Longrightarrow \mathbf{e}(p)/l$, where $p \in M$ and $l$ stands for a unit vector in $L_p$.

Let, as before, $V_i$ denote nonsingular and connected hypersurfaces such that $X = \bigcup V_i$ and let $L_i, \mathbf{e}_i$, $(i = 1, \ldots, k)$ be the associated line bundles and their sections. We define an analytic covering (of degree $2^k$) $\pi \colon \widetilde{M} \to M$ as the fibre product (over $M$) of all such $\pi_{L_i}$. Let $f_i$ be analytic functions on $M$ defining $\widetilde{V}_i = \pi^{-1}(V_i)$ and let $f' = f \circ \pi$. Then:

(1) The generators of the covering group of $\pi$, given by the generators of the covering groups of $\pi_{L_i}$, change the signs of $f_i$;
(2) $f'$ is globally normal crossings and

$$f' = f_1^{a_1} f_2^{a_2} \cdots f_k^{a_k} g,$$

where the positive integers $a_i$ are equal to the generic multiplicities of $f$ along $V_i$ and $g$ is invertible. If, as we have assumed, $f$ is nonnegative, then all $a_i$ are even and $g$ is positive. (In fact, we can get rid of $g$ by replacing $f_k^{a_k}$ by $(f_k g^{1/a_k})^{a_k}$).

**Definition 3.2.** By the *the cut of $M$ along $X$* we mean

$$M' = Closure(\{p \in \widetilde{M};\ f_i(p) > 0, i = 1, \ldots, k\}).$$

We also define

$$\widetilde{X} = \pi^{-1}(X), \quad \widetilde{V}_i = \pi^{-1}(V_i), \quad X' = \widetilde{X} \cap M'.$$

Note that both $M'$, $X'$ are local quadrants so they posses canonical stratifications with normal crossings. Also $X$ has such stratification which can be described as follows. For each locally normal crossings $f\colon M \to \mathbf{R}$ given by (3.1) in a neighbourhood $U$ of $p \in M$ we define

$$k(p) = \#\{a_i : a_i > 0\}.$$

Then the number of local components of $U \setminus X$ at $p$ is exactly $2^{k(p)}$ and $X = \{p \in M; k(p) > 0\}$. The canonical stratification of $X$ is then defined by the sets of the form $\{p \in X; k(p) = const\}$. The stratifcation of $M$ given by the canonical stratification of $X$ and the open stratum $M \setminus X$ we call *associated to X*. We have easily the following lemma.

**Lemma 3.3.** *The induced map $\pi' = \pi_{|M'}\colon M' \to M$ is finite and surjective. The restriction of $\pi'$ to each stratum of the canonical stratification of $M'$ is an analytic covering over a stratum of the associated stratification of $M$. For each $p \in M$*

$$\#\{(\pi')^{-1}(p)\} = 2^{k(p)}.$$

*Since $(\pi')^{-1}(X) = X'$, the same is true for $\pi'|_{X'}$.*

**3.3. Construction of vector fields.** Now we construct an analytic vector field on $\widetilde{M}$ such that the generated flow gives a retraction of $M'$ onto $X'$. This vector field has properties similar to those of $\mathbf{w}$ of (3.2). First we construct some auxiliary vector fields. Let, as above, $a_i$ denote the generic multiplicity of $f$ along $V_i$.

**Lemma 3.4.** *For each $i_0 \in \{1, \ldots, k\}$ there is a vector field $\mathbf{v} = \mathbf{v}_{i_0}$ on $\widetilde{M}$ such that:*

(a) $\qquad\qquad\qquad \partial f_i/\partial \mathbf{v} \equiv 0$ *on* $\widetilde{V}_i$, *for* $i \neq i_0$;

(b) $\qquad\qquad\qquad \partial f_{i_0}/\partial \mathbf{v} \equiv 1/a_{i_0}$ *on* $\widetilde{V}_{i_0}$;

(c) $\qquad\qquad\qquad f_{i_0}\partial f'/\partial \mathbf{v} \equiv f'$ *everywhere on* $M$.

*Proof.* First we show the existence of such $\mathbf{v}$ locally around each $p \in \widetilde{M}$. We consider 3 different cases.

*Case 1.* $p \notin \widetilde{X}$. Then the conditions (a) and (b) are empty. Since $\tilde{f}$ is regular near $p$ we can take

$$\mathbf{v} = \frac{f'}{f_{i_0}\|grad\, f'\|^2}\, grad\, f'$$

*Case 2.* $p \in \widetilde{X} \setminus \widetilde{V}_{i_0}$. Let $k' = k(\pi(p))$ and let $\{f_{i_j} | j = 1, \ldots, k'\}$ be the set of all $f_i$ vanishing at $p$. Take a local system of coordinates $x_1, \ldots, x_n$ around $p$ such that $x_j = f_{i_j}$ for $j = 1, \ldots, k'$. Then near $p$

(3.3)
$$f'(x) = \prod_{j=1}^{k'} x_j^{a_{i_j}} h(x), \text{ where } h \text{ is invertible.}$$

Then we can take

$$\mathbf{v}(x) = \bigl(f_{i_0} \sum_{j=1}^{k'}(a_{i_j} + x_j \partial(lnh)/\partial x_j)\bigr)^{-1} \bigl(\sum_{j=1}^{k'} x_j \partial/\partial x_j\bigr).$$

*Case 3.* $p \in \widetilde{V}_{i_0}$. We take $k'$ and a system of coordinates $x_1, \ldots, x_n$ as in Case 2 and assume that $x_1 = f_{i_0}$. Then again $f'$ is given near $p$ by (3.3) and we may take

$$\mathbf{v}(x) = \bigl(a_{i_0} + x_1 \partial(lnh)/\partial x_1\bigr)^{-1} \partial/\partial x_1.$$

To construct $\mathbf{v}$ globally we apply the method of [Kuo1]. Let $\mathcal{F}$ be a coherent sheaf of analytic vector fields on $\widetilde{M}$ satisfying (a). Let $\mathcal{F}_0$ be a coherent subsheaf of $\mathcal{F}$ of those vector fields which satisfy

(b') $\qquad\qquad\qquad \partial f_{i_0}/\partial \mathbf{v} \equiv 0 \text{ on } \widetilde{V}_{i_0};$

(c') $\qquad\qquad\qquad \partial f'/\partial \mathbf{v} \equiv 0 \text{ everywhere on } \widetilde{M}.$

Then, by Cartan's Theorem B [Car] the following sequence is exact

$$0 \longrightarrow H^0(\mathcal{F}_0) \longrightarrow H^0(\mathcal{F}) \longrightarrow H^0(\mathcal{F}/\mathcal{F}_0) \longrightarrow 0.$$

The local vector fields, we have just constructed, yield a global section of $\mathcal{F}/\mathcal{F}_0$, which by exactness is the image of a global section of $\mathcal{F}$. This is the desired vector field. $\square$

**Lemma 3.5.** *Let $\mathbf{w} = -\sum_{i=1}^{k} \mathbf{v}_i$, where $\mathbf{v}_i$ are vectors fields satisfying the statement of Lemma 3.4. Then there exists a neighbourhood $U$ of $M'$ such that:*

(1) *$\mathbf{w}$ is nonsingular in $U$;*
(2) *on $X'$ the vector field $\mathbf{w}$ is directed outside $M'$, that is, in any local coordinates $x_1, \ldots, x_n$ on $\widetilde{M}$ such that locally $M' = \{x_1 \geq 0, \ldots, x_{k'} \geq 0\}$*

$$\mathbf{w} = \sum_{i=1}^{k'} w^i(x) \partial/\partial x_i + O(\|x\|) \text{ and each } w^i(0) < 0;$$

(3) *for each $p \in U$ there is exactly one $\delta(p) \in \mathbf{R}$ such that if $\gamma(t)$ is an integral curve of $\mathbf{w}$ and $\gamma(0) = p$, then $\gamma(\delta(p)) \in X'$.*

*Proof.* It is easy to see that (2) holds. Therefore $\mathbf{w}$ does not vanish on $X'$. By (c) of Lemma 3.4

$$\partial f'/\partial \mathbf{w} = (\sum_{i=1}^{k} 1/f_i)f', \tag{3.4}$$

and thus $\mathbf{w}$ does not vanish on $M' \setminus X'$. This shows (1). If $\gamma(t)$ is an integral curve of $\mathbf{w}$ and $\gamma(0) \in M'$, then by (3.4) $f'$ decreases along $\gamma$. Therefore, since $f'$ is proper, $\gamma$ has to hit $X'$ and does it only once by (2). This shows (3). □

Let $\Phi(p,t)$ be the flow generated by $\mathbf{w}$. Then $\Phi$ defines a retraction of $r' \colon M' \to X'$ by

$$r'(p) = \Phi(p, (\delta(p))).$$

By (1) of Lemma 3.5 we may assume that $\Phi$ is analytic. Thus the properties of $r'$ depend on the properties of $\delta$.

**Lemma 3.6.** *$\delta(p)$ is a continuous and subanalytic function on $M'$. Moreover, there is a stratification $\mathcal{S}_{M'}$ of $M'$ with normal crossings, which on $X'$ coincides with the canonical stratification of $X'$, and such that $\delta|_{M'}$ is piecewise analytic with respect to $\mathcal{S}_{M'}$.*

*Proof.* We shall substratify $M'$ in a way similar we substratified the quadrant $U_1$ in the local case.

Fix $p \in X'$ and take local coordinates $x_1, \ldots, x_n$ at $p$ satisfying (2) of Lemma 3.5. Let $\delta_i(x)$, $i = 1, \ldots, k'$, be an analytic function such that $\gamma(x, \delta_i(x)) \in \{x_i = 0\}$. Then $\{\delta_i = 0\} = \{x_i = 0\}$ so we can replace the coordinate $x_i$ by $\delta_i(x)$ for every $i = 1, \ldots, k'$. Then, after replacing the remaing coordinate functions $x_j$ by $x_j(\Phi(x, \delta_1(x)))$, $j = k'+1, \ldots, n$,

$$\delta(x) = min\,\{x_i;\ i = 1, \ldots, k'\};$$
$$\mathbf{w} = -(\partial/\partial x_1 + \ldots + \partial/\partial x_{k'}).$$

Therefore, $\mathbf{w}$ looks exactly like in (3.2). As the strata of $\mathcal{S}_{M'}$ we take the inverse images by $r'$ of the strata of the canonical stratification of $X'$. This ends the proof of the lemma. □

The proof of Lemma 3.6 shows that in some local coordinate system $\Phi$ and $r'$ looks exactly the same as in the local case considered in Section 3.1. The strata of $\mathcal{S}_{M'}$ are transverse to nonsigular fibres $M'_c = {f'}^{-1}(c) \cap M'$, $c \in (0, \varepsilon)$, and so induce on $M'_c$ a stratification with normal crossings.

Let

$$R' \colon M' \to X' \times [0, \varepsilon)$$

be given by $R'(p) = (r'(p), f'(p))$. By construction, $R'$ is a homeomorphism and $R'$ trivializes $f'|_{M'}$. In the local coordinates constructed in the proof of Lemma 3.6 this trivialization is given by

$$(3.5) \qquad R'(x) = (x_1 - \delta(x), \ldots, x_{k'} - \delta(x), x_{k'+1}, \ldots, x_n, x_1^{a_1} \cdots x_{k'}^{a_{k'}} g(x)),$$

where $g(x)$ is invertible and $\delta(x) = min\{x_i; i = 1, \ldots, k'\}$. Therefore, by the same argument as in the local case, we obtain the following.

**Proposition 3.7.** *The trivialization $R' \colon M' \to X' \times [0, \varepsilon)$ and the retraction $r'$ satisfy*
  (1) *$R'$ sends the strata of $\mathcal{S}_{M'}$ onto the strata of the product stratification of $X' \times [0, \varepsilon)$ and is piecewise analytic;*
  (2) *$r'$ is a piecewise analytic retraction. Restricted to a nonsingular fibre $M'_c = {f'}^{-1}(c) \cap M'$, $c \in (0, \varepsilon)$, $r'$ is a piecewise analytic isomorphism (which, up to an analytic isomorphism does not depend on the choice of $c$).*

We push so $R'$ and $r'$ onto $M$ and define

$$R \colon M \to X \times [0, \varepsilon), \quad r \colon M \to X,$$

by $r(\pi'(p)) = \pi'(r'(p))$ and $R(\pi'(p)) = \pi'(R'(p)) = (r(\pi'(p)), f(\pi'(p)))$. We define a stratification with normal crossings $\mathcal{S}_M$ of $M$ by taking the images of the strata of $\mathcal{S}_{M'}$. Recall that on $X$, and so on $X \times [0, \varepsilon)$, there is a canonical stratification with normal crossings.

**Proposition 3.8.** *Both $R$ and $r$ are stratified piecewise analytic maps. Moreover, the induced specialization map $r_c = r|_{f^{-1}(c)}$, $c \in (0, \varepsilon)$, is equivalent (up to a piecewise analytic isomorphism) to the projection $\pi'|_{X'} \colon X' \to X$. In particular, $r_c$ is an analytic covering over each stratum with the fibre $r^{-1}(p) \simeq (S^0)^{k(p)}$ (so $\#\{r^{-1}(p)\} = 2^{k(p)}$).*

If $f$ is not everywhere nonnegative, as we assumed above, then we have to modify slightly our construction and divide $M'$ into two disjoint local quadrants

$$M'_+ = Closure(\{p \in \widetilde{M}; g(p) > 0\}), \quad M'_- = Closure(\{p \in M'; g(p) < 0\}).$$

Then we retract independently $M'_+$ onto $X'_+ = X \cap M'_+$ and $M'_-$ onto $X'_- = X \cap M'_-$.

**3.4. General case.** Suppose now that $f: \mathcal{X} \to \mathbf{R}$ is a proper real analytic function with smooth generic fibres. Again, for exposition reason, we assume that $f$ is nonnegative. Let $\mathcal{X}_0 = f^{-1}(0)$ and let $\sigma: M \to \mathcal{X}$ be a modification of $\mathcal{X}$ such that:

(1) $M$ is nonsingular and $\sigma$ induces an isomorphism between $M \setminus \sigma^{-1}(\mathcal{X}_0)$ and $\mathcal{X} \setminus \mathcal{X}_0$;
(2) $f \circ \sigma$ is globally normal crossings.

Let $\pi': M' \to M$ be the cutting along $X = \sigma^{-1}(\mathcal{X}_0)$. Applying the construction of Section 3.3 to $M$ we get the following diagram

$$\begin{array}{c} X' \subset M' \\ \downarrow \quad \downarrow \pi' \\ X \subset M \\ \downarrow \quad \downarrow \sigma \\ \mathcal{X}_0 \subset \mathcal{X} \end{array}$$

and the trivialization $R': M' \to X' \times [0, \varepsilon)$ and retraction $r': M' \to X'$. Denote $\sigma \circ \pi'$ by $\sigma'$. We define the retraction

(3.6) $\qquad r: \mathcal{X} \to \mathcal{X}_0, \text{ by } r(\sigma'(p)) = \sigma'(r'(p)).$

Then Proposition 3.7 implies the following properties of $r$.

**Theorem 3.9.** *Let $f: \mathcal{X} \to \mathbf{R}$ be a proper real analytic function with smooth generic fibres and such that $f$ is nonnegative. Let $\mathcal{X}_0 = f^{-1}(0)$. Fix a modification $\sigma: M \to X$ such that: $M$ is nonsingular, $\sigma$ induces an isomorphism between $M \setminus \sigma^{-1}(\mathcal{X}_0)$ and $\mathcal{X} \setminus \mathcal{X}_0$, and $f \circ \sigma$ is globally normal crossings. Then the retraction $r: \mathcal{X} \to \mathcal{X}_0$ given by (3.6) and the induced specialization $r_c: \mathcal{X}_c \to \mathcal{X}_0$, $c > 0$, satisfy*

(1) *$r$ and $r_c$ are subanalytic (that is their graphs are subanalytic);*
(2) *up to an analytic isomorphism $r_c$ does not depend on $c$.*
(3) *up to a piecewise analytic isomorphism $r_c$ is equivalent to $\sigma'|_{X'}: X' \to \mathcal{X}_0$. In particular, if $\mathcal{S}$ is an analytic stratification of $\mathcal{X}_0$ such that $\sigma'|_{X'}$ is locally topologically trivial over the strata of $\mathcal{S}$, then $r_c$ is also locally topologically trivial over the strata of $\mathcal{S}$.*

*Remark 3.10.* Fix $p \in \mathcal{X}_0$ and a local embedding of $(\mathcal{X}, p)$ in $\mathbf{R}^N$. Then, by the Milnor fibre of $f$ at $p$ we mean

$$F(p) = f^{-1}(c) \cap \mathcal{X} \cap B(p, \varepsilon),$$

where $B(p, \varepsilon)$ is a small ball centered at $p$ with radius $\varepsilon > 0$, and $c$ is chosen such that $0 < \delta \ll \varepsilon$. It is well-known that the homotopy type of $F(p)$ depends neither on the choice of $c, \varepsilon$ nor on the embedding $(\mathcal{X}, p) \subset \mathbf{R}^N$ (if $f$ is not nonnegative then we have to distinguish positive and negative levels of $f$). Let $r_c \colon \mathcal{X}_c \to \mathcal{X}_0$ be the specialization constructed above. Then, for every $p \in \mathcal{X}_0$, $r_c^{-1}(p) \simeq {\sigma'}^{-1}(p)$ are homotopically equivalent $F(p)$

$$(3.7) \qquad r_c^{-1}(p) \simeq {\sigma'}^{-1}(p) \sim F(p).$$

Indeed, to show (3.7) we consider the space $M'$ constructed in Section 3.3 and we follow the notation of that section. Define $\varphi \colon M' \to \mathbf{R}$ by

$$\varphi(x) = \|\sigma'(x) - p\|^2.$$

Note that $\varphi$ is subanalytic (it is even the restriction of an analytic function on $\widetilde{M}$). The sets $M'(p, \varepsilon) = {\sigma'}^{-1}(B(p, \varepsilon)) = \varphi^{-1}([0, \varepsilon))$, for $\varepsilon > 0$ varying, form a system of subanalytic neighbourhoods of $X'(p) = {\sigma'}^{-1}(p)$. For $\varepsilon > 0$ sufficiently small, by standard argument as for instance triangulability, the following inclusions are homotopy equivalences

$$X'(p) \hookrightarrow X' \cap M'(p, \varepsilon) \hookrightarrow M'(p, \varepsilon).$$

Morever, if $\varepsilon > 0$ is small enough, then by Sard's Theorem the levels $\varphi^{-1}(\varepsilon)$ are smooth and transverse to the strata of $\mathcal{S}_{M'}$. Thus, we may modify the argument of the proof of Lemma 3.5 and require that the vector field $\mathbf{w}$ is tangent to one level, say, $\varphi^{-1}(\varepsilon_0)$. Then, the retraction along the integral curves of this new vector field preserves $M'(p, \varepsilon_0)$ and sends ${f'}^{-1}(c) \cap M'(p, \varepsilon_0)$ onto $X' \cap M'(p, \varepsilon_0)$. This shows that the inclusion

$${f'}^{-1}(c) \cap M'(p, \varepsilon_0) \hookrightarrow M'(p, \varepsilon_0)$$

is a homotopy equivalence. So is the inclusion

$${r'_c}^{-1}(X(p)) \hookrightarrow M'(p, \varepsilon_0)$$

since it is homotopic to $X'(p) \hookrightarrow M'(p, \varepsilon_0)$. Hence

$$(3.8) \qquad {r'_c}^{-1}(X(p)) \hookrightarrow {f'}^{-1}(c) \cap M'(p, \varepsilon_0)$$

is also a homotopy equivalence. The push down of (3.8) by $\sigma'$ shows that

$$r_c^{-1}(p) = \sigma'({r'_c}^{-1}(X(p))) \hookrightarrow \sigma'({f'}^{-1}(c) \cap M'(p, \varepsilon_0)) = F(p)$$

is a homotopy equivalence, as claimed.

## 4. Complex case

Let $M$ be a complex connected manifold of dimension $n$ and let

$$f\colon M \to \Delta$$

be a complex analytic proper morphism from $M$ onto a small disc $\Delta \subset \mathbf{C}$ centered at 0. We assume that 0 is the only (possible) critical value of $f$ and that the zero set $X$ of $f$ is normal crossings that is it is the union of nonsingular hypersurfaces intersecting transversally. Then, similarly to the real-analytic case, each point of $p \in M$ admits a coordinate neighbourhood $U$ with coordinates $z = (z_1, \ldots, z_n)$ such that

$$f(z) = z_1^{a_1} z_2^{a_2} \cdots z_n^{a_n} g(z),$$

where $a_i \in \mathbf{N} \cup \{0\}$ and $g$ nowhere vanishes in $U$. Let $k(p) = \#\{a_i : a_i \neq 0\}$.

Our goal is to construct a nice subanalytic retraction $r\colon M \to X$.

**4.1. Basic idea.** We follow the idea of A'Campo [A'C] and take the real-analytic oriented blowing up of $X$ in $M$. Then we apply the real analytic construction of Section 3:

(i) Let $\tau\colon M_{\mathbf{R}} \to M$ be the blowing-up of $M$ along $X$ considered as a real analytic subspace of $M$ (more precisely, we blow up the sheaf of ideals of all real analytic functions vanishing on $X$). Then, the exceptional divisor $X_{\mathbf{R}} = \tau^{-1}(X)$ is the zero set of a real analytic normal crossings $f_{\mathbf{R}} = |f \circ \sigma_M|^2$.

(ii) We apply the real analytic case to $M_{\mathbf{R}}$ and $f_{\mathbf{R}}$. The cutting $M'_{\mathbf{R}}$ of $M_{\mathbf{R}}$ along $X_{\mathbf{R}}$ can be interpreted as 'the oriented blowing-up' of $M$ along $X$ (see, for instance, [A'C]). In particular, let $\tau_\Delta\colon \hat\Delta \to \Delta$ be the real analytic blowing up at the origin and let $\pi_\Delta\colon \Delta' \to \hat\Delta$ be the cutting along the exceptional divisor. Then, $\hat\Delta$ is a Möbious strip and $\Delta'$ is the standard band $\simeq [0, \varepsilon] \times S^1$ with coordinates $\rho, \alpha$ that are nothing else but the polar coordinates $(\rho, \alpha)$ on $\Delta$. As we show below, there are uniquely determined maps $f'$, $\hat f$ making the following diagram commutative.

(4.1)
$$\begin{array}{ccccc}
X'_{\mathbf{R}} \subset M'_{\mathbf{R}} & \xrightarrow{f'} & \Delta' & \xrightarrow{\alpha} & S^1 \\
\downarrow \quad \downarrow{\pi'} & & \downarrow{\pi'_\Delta} & & \\
X_{\mathbf{R}} \subset M_{\mathbf{R}} & \xrightarrow{\hat f} & \hat\Delta & & \\
\downarrow \quad \downarrow{\tau} & & \downarrow{\tau_\Delta} & & \\
X \subset M & \xrightarrow{f} & \Delta & &
\end{array}$$

(iii) In fact, we may apply the construction of Section 3.3 to $f_{\mathbf{R}} = |f \circ \sigma_M|^2$ fiberwise over $S^1$. This is possible since the strata of canonical stratification of $X'_{\mathbf{R}}$ project submersively by $\alpha$ onto $S^1$. Therefore, we can choose the vector field $\mathbf{w}$ of Lemma 3.5 tangent to the levels of $\alpha$. In particular, we may assume that the retraction $r' \colon M'_{\mathbf{R}} \to X'_{\mathbf{R}}$ preserves $\alpha$.

(iv) Then we push $r'$ down to $M$ and define the retraction $r \colon M \to X$. To study the specialization $r_\delta = r|_{f^{-1}(\delta)} \colon X_\delta \to X$, $\delta \in \Delta^* = \Delta \setminus \{0\}$, we first consider $f$ restricted to the levels of $f_{\mathbf{R}}$. Let $X_{c,\mathbf{R}} = f_{\mathbf{R}}^{-1}(c)$ for small real $c > 0$ and let $r_{c,\mathbf{R}} = r|_{X_{c,\mathbf{R}}}$. We have the diagram

(4.2)
$$\begin{array}{c} r_{c,\mathbf{R}} \colon X_{c,\mathbf{R}} \longrightarrow X \\ {\scriptstyle \alpha = arg(f)} \downarrow \\ S^1 \end{array}$$

Thus the complex specialization $r_\delta$ is the restriction of $r_{c,\mathbf{R}}$ to the levels of $\alpha$. By construction, as we will see, the fibres of $r_{c,\mathbf{R}}$ over $p \in X$ are isomorphic (even real-analytically) to $(S^1)^{k(p)}$ and $\alpha$ restricted to such a fibre is isomorphic to

(4.3) $$(S^1)^{k(p)} \ni (\alpha_1, \ldots, \alpha_{k(p)}) \longrightarrow \alpha = \sum a_i \alpha_i \mod 2\pi.$$

We describe below this construction in details.

**4.2. Local Case.** We assume that $f$ is given by

$$f(z_1, \ldots, z_n) = \prod_{i=1}^{k} z_i^{a_i},$$

where $a_i$ are positive integers. Then we define

$$\tau \colon \hat{\Delta}^k \times \Delta^{n-k} \to \Delta^k \times \Delta^{n-k}$$

by taking the product of individual blowings up for each polidisc in $\Delta^k$. One may check that $\tau$ is the blowing up of the ideal of all real analytic functions vanishing on $X = \{0\} \times \Delta^{n-k}$.

To show the existence and uniqueness $\hat{f}$ of (4.1) we use its coordinate free characterization by the universal property of blowing-up (see for instance [H]). Let $\mathcal{I}$ be the ideal of real analytic functions on $\Delta$ vanishing at the origin. The universal property of blowing-up says that the map $\hat{f}$ completing the diagram (4.1) exists if and only if $(f \circ \tau)^* \mathcal{I}$ is an invertible sheaf of ideals. This follows from the lemma below.

**Lemma 4.1.** *Let $\tau: \hat{\Delta}^k \to \Delta^k$ be the product of real-analytic blowings up and let $f(z_1, \ldots, z_k) = \prod_{i=1}^{k} z_i^{a_i}$. Then, the ideal $(Re\, f \circ \tau, Im\, f \circ \tau)$ is invertible in $\hat{\Delta}^k$.*

*Proof.* Let $z_j = x_j + iy_j$. Then, by definition of $\tau$, the ideals $\mathcal{I}_j = (x_j \circ \sigma, y_j \circ \sigma)$ are invertible on $\hat{\Delta}^k$. Take, for instance, such $p \in \hat{\Delta}^k$ that all $\mathcal{I}_j$ are generated by $x_j \circ \sigma$. Then
$$\prod_j (x_j + iy_j)^{a_j} = (\prod_j x_j^{a_j})(\prod_j (1 + iy_j/x_j)^{a_j}) = (\prod_j x_j^{a_j})Z,$$
and either $Re\, Z$ or $Im\, Z$ (or both) is nonzero. For general $p \in \hat{\Delta}^k$ the proof is similar. $\square$

To show the existence of $f': \Delta'^k \times \Delta^{n-k} \to \Delta'$ we just note that $f'$ given by

$$(4.4) \quad f'((\rho_1, \alpha_1), \ldots, (\rho_k, \alpha_k), z_{k+1}, \ldots, z_n) = (\prod_{i=1}^{k} \rho_i^{a_i}, a_1\alpha_1 + \cdots + a_k\alpha_k).$$

closes the diagram (4.1). The uniqueness of $f'$ is obvious.

**4.3. Global Construction.** Assume now that $f$ is global normal crossing. By the existence and uniqueness in the local case (Section 4.2 above) there are uniquely determined maps $\hat{f}$ and $f'$ closing the diagram (4.1). We may apply the real analytic case (Section 3.3 above) and construct the retraction $r': M'_{\mathbf{R}} \to X'_{\mathbf{R}}$. Moreover, by the local description (4.4) of $f'$ it is clear that the strata of canonical stratification of $X'_{\mathbf{R}}$ project submersively (via $\alpha$) on $S^1$. This shows that we may require the vector field $\mathbf{w}$ of Lemma 3.5 to be tangent to the levels of $\alpha$ and the induced by $\mathbf{w}$ retraction $r': M'_{\mathbf{R}} \to X'_{\mathbf{R}}$ to preserve $\alpha$. Again we define $r: M \to X$, by pushing $r'$ down, that is $r(\tau \circ \pi'(p)) = \tau \circ \pi'(r'(p))$. We emphasize the following properties of $r$ which follows directly from the construction.

**Proposition 4.2.** *The retraction $r: M \to X$ and the induced $r_{c,\mathbf{R}}: X_{c,\mathbf{R}} \to X$, $r_\delta: X_\delta \to X$ satisfy:*
  (1) *$r$, $r_{c,\mathbf{R}}$, and $r_\delta$ are subanalytic;*
  (2) *up to an analytic isomorphism $r_{c,\mathbf{R}}$ and $r_\delta$ do not depend on the choice of $c$ and $\delta$ respectively;*
  (3) *up to a piecewise analytic isomorphism over $S^1$, the diagram (4.2) is equivalent to*
$$\begin{array}{c} X'_{\mathbf{R}} \xrightarrow{\tau \circ \pi'} X \\ \alpha \downarrow \phantom{XXX} \\ S^1 \phantom{XXXX} \end{array}$$

**4.4. General case.** Similarly to Section 3.4 we generalize the above construction to the case of $f\colon \mathcal{X} \to \Delta$ being a proper complex analytic function with smooth generic fibres. First we fix a modification $\sigma\colon M \to \mathcal{X}$ of $\mathcal{X}$ satisfying:

(1) $M$ is nonsingular and $\sigma$ induces an isomorphism between $M \setminus \sigma^{-1}(\mathcal{X}_0)$ and $\mathcal{X} \setminus \mathcal{X}_0$;
(2) $f \circ \sigma$ is globally normal crossings.

Apply to $f \circ \sigma\colon M \to \Delta$ the construction of Section 4.1 and consider

$$\begin{array}{ccc} X'_{\mathbf{R}} & \subset & M'_{\mathbf{R}} \\ \downarrow & & \downarrow {\scriptstyle \tau' = \tau \circ \pi'} \\ X & \subset & M \\ \downarrow & & \downarrow {\scriptstyle \sigma} \\ \mathcal{X}_0 & \subset & \mathcal{X} \end{array}$$

In Section 4.1 we constructed $r'\colon M'_{\mathbf{R}} \to X'_{\mathbf{R}}$ which we now push down to a retraction $r\colon \mathcal{X} \to \mathcal{X}_0$. Similarly we define $r_{c,\mathbf{R}}\colon X_{c,\mathbf{R}} \to X$, $r_\delta\colon X_\delta \to X$. The following properties of $r$, $r_{c,\mathbf{R}}$, and $r_\delta$ follows from Proposition 4.2.

**Theorem 4.3.** *Let $f\colon \mathcal{X} \to \Delta$ be a proper complex analytic function with smooth generic fibres and let $\mathcal{X}_0 = f^{-1}(0)$. Fix a modification $\sigma\colon M \to X$ such that: $M$ is nonsingular, $\sigma$ induces an isomorphism between $M \setminus \sigma^{-1}(\mathcal{X}_0)$ and $\mathcal{X} \setminus \mathcal{X}_0$, and $f \circ \sigma$ is globally normal crossings. Let $\tau' = \tau \circ \pi'\colon M'_{\mathbf{R}} \to M$ be the oriented real blowing up of $X = (f \circ \sigma)^{-1}(0)$. Then the constructed above $r$, $r_{c,\mathbf{R}}$, and $r_\delta$ satisfy*

(i) *$r$, $r_{c,\mathbf{R}}$, and $r_\delta$ are subanalytic;*
(ii) *up to an analytic isomorphism $r_{c,\mathbf{R}}$ and $r_\delta$ does not depend on the choice of $c$ and $\delta$ respectively.*
(iii) *up to a piecewise analytic isomorphism over $\mathbf{S}^1$ the diagram*

$$\begin{array}{ccc} r_{c,\mathbf{R}}\colon \mathcal{X}_{c,\mathbf{R}} & \longrightarrow & \mathcal{X}_0 \\ {\scriptstyle arg(f)} \downarrow & & \\ S^1 & & \end{array}$$

*is equivalent to*

$$\begin{array}{ccc} X'_{\mathbf{R}} & \xrightarrow{\sigma' = \sigma \circ \tau'} & \mathcal{X}_0 \\ {\scriptstyle \alpha} \downarrow & & \\ S^1 & & \end{array}$$

*where $\alpha$ is the map constructed in Section 4.1.*

**4.5. Corollaries.**

Fix $p \in \mathcal{X}_0$ and a local embedding of $(\mathcal{X}, p)$ in $\mathbf{C}^N$. We define the Milnor fibration of $f$ at $p$ by

$$f^{-1}(S_c^1) \cap \mathcal{X} \cap B(p, \varepsilon) \xrightarrow{arg(f)} S^1,$$

where $S_c^1 = \{z \in \mathbf{C}; |z| = c\}$, $B(p, \varepsilon)$ is a small ball centered at $p$ with radius $\varepsilon > 0$, and $c$ is a real number such that $0 < c \ll \varepsilon$. The fibre of this fibration is called the Milnor fibre. Applying the argument of Remark 3.10, fibrewise with respect to $\alpha$, we get the following description of Milnor fibration in terms of the above construction.

**Corollary 4.4.** *Up to a homotopy equivalence the Milnor fibration of $f$ at $p$ is described by*

$$r_c^{-1}(p) \xrightarrow{arg(f)} S^1,$$

*or equivalently by*

$$(\sigma \circ \tau')^{-1}(p) \xrightarrow{\alpha} S^1.$$

Now we show that the constructed above retraction satisfies the properties requested by Deligne in [De].

Let $\tilde{\Delta} \xrightarrow{s} \Delta^*$ be the universal covering of $\Delta^* = \Delta \setminus \{0\}$ with the covering group $G = \pi_1(\Delta^*)$. Then $s$ induces a covering $\tilde{s}$ closing the diagram

(4.5)
$$\begin{array}{ccc} \widetilde{\mathcal{X} \setminus \mathcal{X}_0} & \xrightarrow{\tilde{f}} & \tilde{\Delta} \\ \downarrow{\tilde{s}} & & \downarrow{s} \\ \mathcal{X} \setminus \mathcal{X}_0 & \xrightarrow{f} & \Delta^* \end{array}$$

Since $\tilde{\Delta}$ is contractible $\tilde{f}$ is topologically trivial. Let $\Gamma: \widetilde{\mathcal{X} \setminus \mathcal{X}_0} \to \tilde{\Delta} \times \mathcal{X}_\delta$ be a $G$-equivariant trivialization of $\tilde{f}$. Following [De] we say that the retraction $r: \mathcal{X} \to \mathcal{X}_0$ and the trivialization $\Gamma$ *are compatible* if $r \circ \tilde{s} \circ \Gamma^{-1}: \tilde{\Delta} \times X_\delta \to \mathcal{X}_0$ coincides with $r_\delta$.

**Corollary 4.5.** *There exists a $\pi_1(\Delta^*)$-equivariant trivialization $\Gamma$ of $\tilde{f}$ which is compatible with retraction $r$ constructed above.*

*Proof.* Let $\hat{s}: \mathbf{R} \to S^1$ be the universal covering and let $\tilde{\alpha}: \widetilde{X'_\mathbf{R}} \to \mathbf{R}$ be induced by

$$\begin{array}{ccc} \widetilde{X'_\mathbf{R}} & \xrightarrow{\tilde{\alpha}} & \mathbf{R} \\ \downarrow & & \downarrow{\hat{s}} \\ X'_\mathbf{R} & \xrightarrow{\alpha} & S^1 \end{array}$$

Recall from Section 3.3 that $M'_{\mathbf{R}} \simeq X'_{\mathbf{R}} \times [0, \varepsilon)$ and the retraction $r'$ corresponds to the projection on the first factor. Let $\rho \colon M'_{\mathbf{R}} \to [0, \varepsilon)$ correspond to the projection on the second factor. Then, via the isomorphism $M'_{\mathbf{R}} \setminus X'_{\mathbf{R}} \xrightarrow{\sigma \circ \tau'} \mathcal{X} \setminus \mathcal{X}_0$,

$$(\rho, \alpha) \colon M'_{\mathbf{R}} \setminus X'_{\mathbf{R}} \to (0, \varepsilon) \times S^1$$

corresponds to $f \colon \mathcal{X} \setminus \mathcal{X}_0 \to \Delta^*$, where we identify $(\rho, \alpha)$ with polar coordinates on $\Delta^*$. Thus, in order to trivialize $\widetilde{\mathcal{X} \setminus \mathcal{X}_0} \xrightarrow{\tilde{f}} \tilde{\Delta}$, it suffices to trivialize $(\rho, \tilde{\alpha})$ on $\widetilde{M'_{\mathbf{R}} \setminus X'_{\mathbf{R}}} = X'_{\mathbf{R}} \times (0, \varepsilon) \times \mathbf{R}$. It is clear that such $\pi_1(S^1)$-equivariant trivialization compatible with the projection onto $X'_{\mathbf{R}}$ can be found. The push down of this trivialization to $\mathcal{X}$ is a trivialization satisfying the required properties. $\square$

As we mentioned in the introduction we would like the retraction $r$ to be locally trivial along the strata of a complex analytic stratification of $\mathcal{X}_0$. By Theorem 4.3 (iii) it suffices that the diagram

(4.6)
$$\begin{array}{c} X'_{\mathbf{R}} \xrightarrow{\sigma' = \sigma \circ \tau'} \mathcal{X}_0 \\ \alpha \downarrow \\ S^1 \end{array}$$

is locally trivial along each stratum. It is not clear to the author that such a complex analytic stratification exists in general, even if we assume that the modification $\sigma \colon M \to \mathcal{X}$ satisfies the condition (1) and (2) of Section 4.4. But this is the case if we make an additional assumption about $\sigma$.

Let $\sigma \colon M \to \mathcal{X}$ be a modification satisfying conditions (1) and (2) of Section 4.4. Then $X = f^{-1}(0)$, being a divisor with normal crossings, has a canonical complex analytic stratification $\mathcal{S}_X$. We assume that, besides (1) and (2), $\sigma$ satisfies the following condition:

(3) *there is a complex analytic stratification $\mathcal{S}$ of $\mathcal{X}_0$ such that the inverse image of the closure of every stratum of $\mathcal{S}$ is the union of some components of $X$ and such that $\sigma|_X \colon X \to \mathcal{X}_0$ sends strata of $\mathcal{S}_X$ submersively onto the strata of $\mathcal{S}$.*

Let $p \in X$ and let $S$ be a stratum of $\mathcal{S}$ such that $\sigma(p) \in S$. Then, by condition (3), $p$ admits a coordinate neighbourhood $U$ with coordinates $z = (z_1, \ldots, z_n)$ such that

$$\sigma^{-1}(S) = \{z_1 z_2 \cdots z_l = 0\}$$

and

$$f(z) = z_1^{a_1} z_2^{a_2} \cdots z_k^{a_k},$$

where $a_i$ are positive integers and $l \leq k$, and the projection $\sigma_S \colon \sigma^{-1}(S) \to S$ is given by
$$\sigma_S(z_1, \ldots, z_n) = (z_{i_1}, \ldots, z_{i_m}),$$
where $\{i_1, \ldots, i_m\} \subset \{k+1, \ldots, n\}$.

Let $U'_{\mathbf{R}} = {\tau'}^{-1}(U)$ and let
$$((\rho_1, \alpha_1), \ldots, (\rho_k, \alpha_k), z_{k+1}, \ldots, z_n)$$
be the coordinates on $U'_{\mathbf{R}}$ (as in Section 4.2). Thus the projection $\sigma'_S \colon (\sigma \circ \tau')^{-1}(S) \to S$ is near $p$ given by
$$((\rho_1, \alpha_1), \ldots, (\rho_k, \alpha_k), z_{k+1}, \ldots, z_n) \to (z_{i_1}, \ldots, z_{i_m}).$$

¿From this local description and (4.3) is clear that $\sigma'_S$ restricted to the levels of $\alpha$ is a submersion. This shows the following corollary.

**Corollary 4.6.** *If the modification $\sigma \colon M \to \mathcal{X}$ satisfies conditions (1), (2), and (3) then the retraction $r$ given by the construction of this section is locally trivial over the strata of the stratification $\mathcal{S}$ of $\mathcal{X}_0$.*

At the end we show that a modification satisfying (3) always exists. Indeed, such a stratification may be constructed by an inductive procedure. First, let $\sigma_1 \colon M_1 \to \mathcal{X}$ be a modification of $\mathcal{X}$ satisfying condition (1) and (2). Let $\operatorname{Reg} \mathcal{X}_0$, $\operatorname{Sing} \mathcal{X}_0$ denote the set of regular and singular points of $\mathcal{X}_0$ respectively. Then there exists a proper analytic subset $Y_1 \subset \mathcal{X}_0$ such that $\operatorname{Sing} \mathcal{X}_0 \subset Y_1$ and the strata of canonical stratification of $X_1 = (f \circ \sigma_1)^{-1}(0)$ project submersively onto $\mathcal{X}_0 \setminus Y_1$. Applying desingularization theorem of Hironaka we may blow up $\sigma_1^{-1}(Y_1)$ so that the resulting modification $M_2 \xrightarrow{\sigma_2} M_1 \xrightarrow{\sigma_1} \mathcal{X}$ satisfies conditions (1) and (2) and the inverse image of $Y_1$ is the union of some components of the global normal crossing $X_2 = (f \circ \sigma_1 \circ \sigma_2)^{-1}(0)$. Then there exists a proper analytic subset $Y_2 \subset Y_1$ such that $\operatorname{Sing} Y_1 \subset Y_2$ and the strata of canonical stratification of $(\sigma_1 \circ \sigma_2)^{-1}(Y_1)$ project submersively onto $Y_1 \setminus Y_2$. It is clear that this procedure leads to a construction of stratification of $\mathcal{X}_0$ satisfying the required properties.

## References

[A'C]  N. A'Campo, *La fonction zêta d'une monodromie*, Comment. Math. Helvetici **50** (1975), 233–248.

[Car]  H. Cartan, *Variétés analytiques réelles et variétés analytiques complexes*, Bull. Soc. Math. France **85** (1957), 77–99.

[Cl]  C. H. Clemens, *Degeneration of Kähler manifolds*, Duke Math. Journ. **44** (1977), no. 2, 215–290.

[CR]   M. Coste, M. Reguiat, *Trivialités en famille*, Real Algebraic Geometry, Lecture Notes in Math., vol. 1524, Springer-Verlag, 1992, pp. 193–204.
[De]   P. Deligne, *La formalisme des cycles évanescents*, Groupes de Monodromie en Géométrie Algébrique (SGA 7 II), Lect. Notes in Math., vol. 340, Springer, Berlin Heidelberg New York, 1973, pp. 82–116.
[GM]   M. Goresky, R. MacPherson, *Morse Theory and Intersection Homology*, Analyse et Topologie sur les Espaces Singuliers, Astérisque **101** (1983), 135–192.
[Gr]   H. Grauert, *On Levi's problem and the imbedding of real-analytic manifolds*, Ann. Math. **68** (1958), no. 2, 460–472.
[H]    H. Hironaka, *Resolution of singularities of an algebraic variety over a field of characteristic zero* I,I*I*, Ann. of Math. **79** (1964), no. 2, 109-326.
[Kuo1] T.-C. Kuo, *On classification of real singularities*, Invent. Math., vol. 82, 1985, pp. 257–262.
[Kuo2] \_\_\_\_\_\_, *A natural equivalence relation on singularities*, Singularities, Banach Center Publications (S. Łojasiewicz, ed.), vol. XX, PWN, Warszawa, 1988, pp. 239–243.
[Sh1]  M. Shiota, *Piecewise linearization of real-valued subanalytic functions*, Trans. Amer. Math. Soc. **312** (2), 663–679.
[Sh2]  \_\_\_\_\_\_, *Piecewise linearization of subanalytic functions II*, Real analytic and algebraic geometry, Lect. Notes Math, vol. 1420, Springer-Verlag, pp. 247–307.

(permanent address) Département de Mathématiques, Université d'Angers, 2 bd. Lavoisier, 49045 Angers, France

*E-mail address*: parus@tonton.univ-angers.fr